\documentclass[11pt,dvips]{article}
\usepackage{amsmath,amsfonts,amssymb,amsthm,amscd,graphicx,verbatim,subfig}
\usepackage[all]{xy}
\usepackage[a4paper,left=25mm]{geometry}
\ifx\pdftexversion\undefined

\usepackage[a4paper,colorlinks,link
=black,filecolor=black,citecolor=black,urlcolor=black,pdfstartview=FitH]{hyperref}
\else

\usepackage[a4paper,colorlinks,linkcolor=black,filecolor=black,citecolor=black,urlcolor=black,pdfstartview=FitH]{hyperref}
\fi

\font\sixbb=msbm6
\font\eightbb=msbm8
\font\twelvebb=msbm10 scaled 1095
\newfam\bbfam
\textfont\bbfam=\twelvebb \scriptfont\bbfam=\eightbb
                           \scriptscriptfont\bbfam=\sixbb
\def\bb{\fam\bbfam\twelvebb}
\newcommand{\Rea}{{\bb R}}

\newcommand{\Int}{{\bb Z}}
\newcommand{\Rat}{{\bb Q}}

\newcommand{\KK}{{\bb K}}

\newtheorem{theorem}{\bf Theorem}
\newtheorem{claim}[theorem]{\bf Claim}
\newtheorem{conjecture}[theorem]{\bf Conjecture}
\newtheorem{proposition}[theorem]{\bf Proposition}
\newtheorem{corollary}[theorem]{\bf Corollary}
\newcommand{\enp}{\begin{flushright} $\Box$ \end{flushright}}
\newcommand{\beq}[0]{\begin{equation}}
\newcommand{\enq}[0]{\end{equation}}

\newcommand{\dn}{\Delta_{n-1}}

\newcommand{\thh}{\tilde{H}}
\newcommand{\supp}{{\rm supp}}
\newcommand{\cf}{{\cal F}}

\newcommand{\tb}{\tilde{\beta}}

\newcommand{\lk}{{\rm lk}}
\newcommand{\str}{{\rm st}}
\newcommand{\costr}{{\rm cost}}

\title{Betti Numbers of Complexes \\
with Highly Connected Links}
\begin{document}
\author{Amir Abu-Fraiha\thanks{Department of Mathematics, Technion, Haifa
    32000, Israel. e-mail: aameer2@gmail.com~.} \and Roy Meshulam\thanks{Department of Mathematics, Technion, Haifa
    32000, Israel. e-mail: meshulam@math.technion.ac.il~. Supported by
     ISF and GIF grants.}}

\maketitle
\pagestyle{plain}
\begin{abstract}
Let $\dn^{(k)}$ denote the $k$-dimensional skeleton of the $(n-1)$-simplex $\dn$ and consider a complex
$\dn^{(k-1)} \subset X \subset \dn^{(k)}$. Let $\KK$ be a field and let $0 \leq \ell <k$. It is shown that if $\thh_{k-\ell-2}(\lk(X,\tau);\KK)=0$ for all $\ell$-dimensional faces $\tau$ of $X$ then
$$
\dim \thh_{k-1}(X;\KK) \leq \frac{\binom{n-1}{\ell} \binom{n-\ell-2}{k-\ell}}{\binom{k+1}{\ell+1}}
$$
with equality iff $\lk(X,\tau)$ is a $(k-\ell-1)$-hypertree for all $\ell$-dimensional simplices $\tau$ of $\dn$.
Examples based on sum complexes show that the bound is asymptotically tight for all fixed $k,\ell$ as $n \rightarrow \infty$.
\end{abstract}

\section{Introduction}

\ \ \
Let $X$ be a simplicial complex on the vertex set $V$. Numerous problems in topological combinatorics ask for estimates on some global invariants of $X$, e.g. its connectivity or Betti numbers, given that $X$ satisfies certain local properties. One remarkable local to global result of this nature is Garland's theorem \cite{Gar73}.
We first recall some definitions. The induced subcomplex of $X$ on $V' \subset V$ is
$X[V']=\{\sigma \in X: \sigma \subset V'\}$.
Denote the star, link and costar of a simplex $\tau \in X$ by
\begin{equation*}
\begin{split}
\str(X,\tau)&=\{\sigma \in X: \sigma \cup \tau \in X\} \\
\lk(X,\tau)&=\{\sigma \in \str(X,\tau): \sigma \cap \tau=\emptyset\} \\
\costr(X,\tau)&=\{\sigma \in X : \sigma \not\supset \tau \}.
\end{split}
\end{equation*}
Let $X^{(j)}$ denote the $j$-th skeleton of $X$ and let $X(j)$ be the family of $j$-dimensional simplices of $X$. Denote $f_j(X)=|X(j)|$.
Assume that $X$ is a pure $k$-dimensional complex, and define a positive weight function on its simplices by
$$c(\sigma)=(k-\dim \sigma)! |\{\tau \in X(k): \tau \supset \sigma\}|.$$
For $\tau \in X$ let $c_{\tau}$ be the induced weight function on $\lk(X,\tau)$ given by $c_{\tau}(\alpha)=c(\tau \cup \alpha)$.
For $-1 \leq j \leq k$ let $C^j(X;\Rea)$ denote the space of real valued $j$-cochains of $X$ and let $d_j:C^j(X;\Rea) \rightarrow C^{j+1}(X;\Rea)$ denote the $j$-th coboundary map
of $X$. Let $d_j^*:C^{j+1}(X;\Rea) \rightarrow C^j(X;\Rea)$ be the adjoint of $d_j$ with respect to the weight function $c$.
Let $L_j=d_{j-1}d_{j-1}^*+d_j^*d_j:C^j(X;\Rea) \rightarrow C^j(X;\Rea)$ be the $j$-th Laplacian of $X$ and let $\mu_j(X)$ denote its minimal eigenvalue.
\begin{theorem}[Garland \cite{Gar73}]
\label{garland}
Let $-1 \leq \ell < k-1$.  If
$\mu_{k-\ell-2}(\lk(X,\tau))>\frac{\ell+1}{k}$ for all $\tau \in X(\ell)$, then $\thh_{k-1}(X;\Rea)=0.$
\end{theorem}
\noindent
Garland's Theorem and its variants have applications in a wide range of areas including representation theory, geometric group theory, hypergraph matching and random complexes
(see e.g. \cite{Gar73,BS97,ABM05,Kahle14}). Here we study the following question that naturally arises in connection with Theorem \ref{garland}:  What can be said  concerning $\thh_{k-1}(X)$ if, instead of $\mu_{k-\ell-2}(\lk(X,\tau))>\frac{\ell+1}{k}$, it is only assumed that $\thh_{k-\ell-2}(\lk(X,\tau))=0$ for all $\tau \in X(\ell)$?
\ \\ \\
Let $\KK$ be a fixed field. Let $\thh_j(X)=\thh_j(X;\KK)$ and  $\tb_j(X)=\dim_{\KK} \thh_j(X;\KK)$ be the reduced homology groups and reduced Betti numbers of $X$ over $\KK$.
For $0 \leq  \ell <k$ and $-1 \leq j$ let
$$
\lambda_{\ell,j}(X)=\sum_{\tau \in X(\ell)} \tb_j(\lk(X,\tau)).
$$
Let
$$
B_{n,k,\ell}=\frac{\binom{n-1}{\ell} \binom{n-\ell-2}{k-\ell}}{\binom{k+1}{\ell+1}}
$$
and
$$F_{n,k,\ell}=\binom{n-1}{k}-B_{n,k,\ell}=\frac{\binom{n}{\ell+1}\binom{n-\ell-2}{k-\ell-1}}{\binom{k+1}{\ell+1}}.
$$
\noindent
Let $\Delta(V)$ denote the simplex on the vertex set $V$. Let $[n]=\{0,\ldots,n-1\}$ and let $\dn=\Delta([n])$ be the $(n-1)$-simplex on $[n]$.
\begin{theorem}
\label{zlink}
If $\dn^{(k-1)} \subset X \subset \dn$ and $0 \leq \ell <k$ then
\begin{equation}
\label{ulink}
\tb_{k-1}(X) \leq \frac{1}{\binom{k+1}{\ell+1}}\lambda_{\ell,k-\ell-2}(X) +B_{n,k,\ell}.
\end{equation}
\end{theorem}
\noindent
Theorem \ref{zlink} implies the following
\begin{corollary}
\label{zlink1}
Suppose $\dn^{(k-1)} \subset X \subset \dn$ satisfies $\lambda_{\ell,k-\ell-2}(X)=0$. Then
\begin{equation}
\label{bzlink}
\tb_{k-1}(X) \leq B_{n,k,\ell}.
\end{equation}
\end{corollary}
\noindent
The next three results concern some aspects of the equality cases in Corollary \ref{zlink1}. A complex $\Delta(V)^{(r-1)} \subset Y \subset \Delta(V)^{(r)}$ is an
{\it $r$-hypertree over $\KK$} on the vertex set $V$ (abbreviated $r$-hypertree for a fixed field $\KK$) if
$\thh_*(Y;\KK)=0$. It is easy to check that $Y$ is an $r$-hypertree iff $f_{r}(Y)=\binom{|V|-1}{r}$ and either $\thh_{r-1}(Y;\KK)=0$ or $\thh_{r}(Y;\KK)=0$.
See Kalai's paper \cite{K83} for further discussion, including a Cayley type formula for the weighted enumeration of rational hypertrees.
\begin{theorem}
\label{zlink2}
Let $0 \leq \ell <k$ and suppose $\dn^{(k-1)} \subset X \subset \dn^{(k)}$ satisfies
$\lambda_{\ell,k-\ell-2}(X)=0.$
Then the following three conditions are equivalent.
\begin{enumerate}
\item[(a)]
$
\tb_{k-1}(X) = B_{n,k,\ell}.
$
\item[(b)]
$\tb_k(X)=0$ and $f_k(X)=F_{n,k,\ell}$.

\item[(c)]
$\lk(X,\tau)$ is a $(k-\ell-1)$-hypertree on $[n]-\tau$ for all $\tau \in \dn(\ell)$.
\end{enumerate}
\end{theorem}
\noindent
The next result asserts that the bound \eqref{bzlink} is asymptotically tight for fixed
$k,\ell$ and $n \rightarrow \infty$.
\begin{theorem}
\label{exacom}
Let $0 \leq \ell < k$ be fixed. Then for any prime number $n > k$ there exists a complex
$\dn^{(k-1)} \subset X_{n,k,\ell} \subset \dn^{(k)}$ such that $\lambda_{\ell,k-\ell-2}(X_{n,k,\ell})=0$ and
$$\tb_{k-1}(X_{n,k,\ell}) \geq (1-O(n^{-1}))B_{n,k,\ell}.
$$
\end{theorem}
\noindent
Finally, we give examples that show the optimality of \eqref{bzlink} for $\ell=0$ and $k \leq 3$.
\begin{theorem}
\label{exactc}
Let $1 \leq k \leq 3$. Then for infinitely many $n$'s there exist complexes
$\dn^{(k-1)} \subset J_{n,k} \subset \dn^{(k)}$ such that $\tb_{k-2}(\lk(J_{n,k},v))=0$ for all $v \in \dn(0)$ and
$$\tb_{k-1}(J_{n,k})= B_{n,k,0}=\frac{1}{k+1} \binom{n-2}{k}.$$
\end{theorem}

The paper is organized as follows: In Section \ref{s:upper} we prove a monotonicity result (Proposition \ref{mono}) that directly implies Theorem \ref{zlink}. The characterization of equality cases (Theorem \ref{zlink2}) is established in Section \ref{s:lkzero}.
In Section \ref{s:sum} we recall the notion of sum complexes and prove an upper bound on the Betti number of their links (Proposition \ref{linkb}). This result is the main ingredient in the proof of Theorem \ref{exacom} given in Section \ref{s:lower}. In Section \ref{s:twoc} we describe the constructions that yield Theorem \ref{exactc}. We conclude in Section \ref{s:conc} with some remarks and open problems.

\section{The Upper Bound}
\label{s:upper}

The main ingredient in the proof of Theorem \ref{zlink} is the following monotonicity result.
\begin{proposition}
\label{mono}
Let $0 \leq \ell <k$. If $X \subset \dn^{(k)}$ and $\sigma \in X(k)$ then
\begin{equation}
\label{xms}
\begin{split}
&\lambda_{\ell,k-\ell-2}(X-\sigma)-\lambda_{\ell,k-\ell-2}(X) \\
&\leq \binom{k+1}{\ell+1}\left(\tb_{k-1}(X-\sigma)-\tb_{k-1}(X) \right).
\end{split}
\end{equation}
\end{proposition}
\noindent
{\bf Proof:} First note that
$$
\tb_{k-1}(X) \leq \tb_{k-1}(X-\sigma) \leq \tb_{k-1}(X)+1.
$$
Let $\tau \in X(\ell)$. If $\tau \subset \sigma$ then
\begin{equation}
\label{lamk1}
\tb_{k-\ell-2}(\lk(X,\tau)) \leq \tb_{k-\ell-2}(\lk(X-\sigma,\tau)) \leq \tb_{k-\ell-2}(\lk(X,\tau))+1.
\end{equation}
On the other hand, if $\tau \not\subset \sigma$ then $\lk(X,\tau)=\lk(X-\sigma,\tau)$.
Summing \eqref{lamk1} over all $\tau \in X(\ell)$ we obtain
\begin{equation}
\label{boundg1}
\lambda_{\ell,k-\ell-2}(X-\sigma) \leq \lambda_{k-\ell-2}(X)+\binom{k+1}{\ell+1}
\end{equation}
Consider two cases:
\\
(i) $\tb_{k-1}(X-\sigma)=\tb_{k-1}(X)+1$. Then by \eqref{boundg1}
\begin{equation*}
\begin{split}
\binom{k+1}{\ell+1}&\left(\tb_{k-1}(X-\sigma)-\tb_{k-1}(X) \right)\\
&=\binom{k+1}{\ell+1}\geq \lambda_{\ell,k-\ell-2}(X-\sigma)-\lambda_{k-\ell-2}(X).
\end{split}
\end{equation*}
Thus \eqref{xms} holds.
\\
(ii) $\tb_{k-1}(X-\sigma)=\tb_{k-1}(X)$. To establish \eqref{xms} it suffices to show that $\lambda_{\ell,k-\ell-2}(X-\sigma)= \lambda_{\ell,k-\ell-2}(X)$, or equivalently
that if $\tau \in \sigma(\ell)$ then
\begin{equation}
\label{equlg}
\tb_{k-\ell-2}(\lk(X,\tau))=\tb_{k-\ell-2}(\lk(X-\sigma,\tau)).
\end{equation}
Consider the decompositions
$$X-\sigma=\costr(X,\tau) \cup \str(X-\sigma,\tau)$$
and
$$X=\costr(X,\tau) \cup \str(X,\tau).$$
Then
$$\costr(X,\tau) \cap \str(X-\sigma,\tau) = \partial\tau * \lk(X-\sigma,\tau)$$
and
$$\costr(X,\tau) \cap \str(X,\tau) = \partial\tau * \lk(X,\tau).$$
Note that $\thh_{k-2}(\partial\tau *Y)\cong \thh_{k-\ell-2}(Y)$ for any $Y$.
Hence by Mayer-Vietoris we obtain a commutative diagram
\begin{scriptsize}
$$
\minCDarrowwidth24pt
\begin{CD}
\thh_{k-1}(\costr(X,\tau)) @>>>
\thh_{k-1}(X-\sigma) @>>>
\thh_{k-\ell-2}(\lk(X-\sigma,\tau)) @>>> \thh_{k-2}(\costr(X,\tau)) @>>> \thh_{k-2}(X-\sigma)  \\
@VV{(i_1)_*}V @VV{(i_2)_*}V @VV{(i_3)_*}V @VV{(i_4)_*}V @VV{(i_5)_*}V \\
\thh_{k-1}(\costr(X,\tau)) @>>>
\thh_{k-1}(X) @>>>
\thh_{k-\ell-2}(\lk(X,\tau)) @>>> \thh_{k-2}(\costr(X,\tau)) @>>> \thh_{k-2}(X)
\end{CD}
$$
\end{scriptsize}
where the rows are exact and the $i_*$'s are induced by inclusion maps. Clearly $(i_1)_*$ and $(i_4)_*$ are the identity maps. As the removal of the $k$-dimensional simplex $\sigma$ does not effect the $(k-2)$-homology, it follows that $(i_5)_*$ is an isomorphism. The assumption $\tb_{k-1}(X-\sigma)=\tb_{k-1}(X)$ implies that
$(i_2)_*$ is an isomorphism. It follows by the $5$-lemma that $(i_3)_*$ is an isomorphism as well, and thus \eqref{equlg} holds.
This completes the proof of \eqref{xms}.
{\enp}
\noindent
{\bf Proof of Theorem \ref{zlink}:}
First note that if $\tau \in \dn(\ell)$ and $\ell<k$ then $\lk(\dn^{(k-1)},\tau)\cong \Delta_{n-\ell-2}^{(k-\ell-2)}$ and that $\tb_j(\Delta_m^{(j)})=\binom{m}{j+1}$.
Secondly, as both $\thh_{k-1}(X)$ and $\thh_{k-\ell-2}(\lk(X,\tau))$ for $\tau \in \dn(\ell)$ depend only on the $k$-dimensional skeleton of $X$, we may assume that
$X \subset \dn^{(k)}$. By repeatedly removing $k$-simplicies from $X$ and using \eqref{xms} it follows that
\begin{equation*}
\label{ddow}
\begin{split}
&\binom{k+1}{\ell+1}\tb_{k-1}(X)-\lambda_{\ell,k-\ell-2}(X) \\
& \leq \binom{k+1}{\ell+1}\tb_{k-1}(\dn^{(k-1)})-\lambda_{\ell,k-\ell-2}(\dn^{(k-1)}) \\
&=\binom{k+1}{\ell+1} \binom{n-1}{k}-\binom{n}{\ell+1}\binom{n-\ell-2}{k-\ell-1} \\
&=\binom{n-1}{\ell} \binom{n-\ell-2}{k-\ell}.
\end{split}
\end{equation*}
{\enp}
\noindent
Theorem \ref{zlink} can be also formulated as the following upper bound on $\tb_k(X)$.
\begin{theorem}
\label{betak}
Let $X \subset \dn^{(k)}$. Then for any $-1 \leq \ell < k$
\begin{equation}
\label{bbkmo}
\binom{k+1}{\ell+1} \tb_k(X) \leq \sum_{\tau \in X(\ell)} \tb_{k-\ell-1}(\lk(X,\tau)).
\end{equation}
\end{theorem}
\noindent
{\bf Proof:}
As both sides of \eqref{bbkmo} do not depend on the $(k-1)$-skeleton of $X$, we may assume that $X  \supset \dn^{(k-1)}$.
The exact sequence for the pair $(X,\dn^{(k-1)})$
$$
0 \rightarrow \thh_k(X) \rightarrow H_k(X,\dn^{(k-1)}) \rightarrow \thh_{k-1}(\dn^{(k-1)}) \rightarrow \thh_{k-1}(X) \rightarrow 0
$$
implies that
\begin{equation}
\label{lexact}
\tb_{k}(X)=\tb_{k-1}(X)+f_k(X)-\binom{n-1}{k}.
\end{equation}
Similarly, for each $\tau \in \dn(\ell)$
\begin{equation}
\label{lexact1}
\tb_{k-\ell-1}(\lk(X,\tau))=\tb_{k-\ell-2}(\lk(X,\tau))+f_{k-\ell-1}(\lk(X,\tau))-\binom{n-\ell-2}{k-\ell-1}.
\end{equation}
Summing \eqref{lexact1} over all $\tau \in \dn(\ell)$ we obtain
\begin{equation}
\label{e:lamdas}
\begin{split}
\lambda_{\ell,k-\ell-1}(X)&=\lambda_{\ell,k-\ell-2}(X)+\sum_{\tau \in \dn(\ell)} f_{k-\ell-1}(\lk(X,\tau)) - \binom{n}{\ell+1}\binom{n-\ell-2}{k-\ell-1} \\
&=\lambda_{\ell,k-\ell-2}(X)+\binom{k+1}{\ell+1}\big(f_k(X)-F_{n,k,\ell}\big).
\end{split}
\end{equation}
Combining \eqref{lexact}, \eqref{ulink} and \eqref{e:lamdas} it follows that
\begin{equation*}
\label{e:equivf}
\begin{split}
\tb_k(X)&=\tb_{k-1}(X)+f_k(X)-\binom{n-1}{k} \\
&\leq \frac{1}{\binom{k+1}{\ell+1}} \lambda_{\ell,k-\ell-2}(X)+B_{n,k,\ell}+f_k(X)-\binom{n-1}{k} \\
&= \frac{1}{\binom{k+1}{\ell+1}} \big(\lambda_{\ell,k-\ell-1}(X)-\binom{k+1}{\ell+1}(f_k(X)-F_{n,k,\ell})\big)+B_{n,k,\ell}+f_k(X)-\binom{n-1}{k} \\
&= \frac{1}{\binom{k+1}{\ell+1}} \lambda_{\ell,k-\ell-1}(X).
\end{split}
\end{equation*}
{\enp}

\section{Characterizations of Equality}
\label{s:lkzero}
In this section we prove Theorem \ref{zlink2}. Denote the support of a $k$-chain $z=\sum_{\sigma \in Y(k)} a_{\sigma} \sigma$ of a complex $Y$ by $\supp(z)=\{\sigma \in Y(k): a_{\sigma} \neq 0\}$.
We shall need the following observation.
\begin{claim}
\label{ktol}
Let $Y \subset \dn^{(k)}$ and let $0 \neq z \in \thh_k(Y)$. If $\sigma \in \supp(z)$ and $\tau \in \sigma(\ell)$
then $\tb_{k-\ell-1}(\lk(Y,\tau)) >0$.
\end{claim}
\noindent
{\bf Proof:} The assumptions imply that $\tb_k(\costr(Y,\tau))<\tb_k(Y)$.
Using, as in the proof of Proposition \ref{mono}, the exact sequence
$$ \thh_k(\costr(Y,\tau)) \rightarrow \thh_k(Y) \rightarrow \thh_{k-\ell-1}(\lk(Y,\tau))$$
it follows that
$$\tb_{k-\ell-1}(\lk(Y,\tau)) \geq \tb_k(Y)-\tb_k(\costr(Y,\tau))>0.$$
{\enp}
\noindent
{\bf Proof of Theorem \ref{zlink2}:}
(a) $\Rightarrow$ (b): We first show that (a) implies $\tb_k(X)=0$. Otherwise choose an inclusion-wise minimal
$\dn^{(k-1)} \subset X \subset \dn^{(k)}$ such that both
$\lambda_{\ell,k-\ell-2}(X)=0$ and $\tb_{k-1}(X)=B_{n,k,\ell}$, but $\tb_k(X)>0$.
Let $0 \neq z \in \thh_k(X)$ and let $\sigma \in \supp(z)$. Then
$\tb_k(X-\sigma)=\tb_k(X)-1$ and hence by \eqref{lexact}
\begin{equation}
\label{bkmo}
\tb_{k-1}(X-\sigma)=\tb_{k-1}(X)=B_{n,k,\ell}.
\end{equation}
\noindent
Proposition \ref{mono} thus implies
\begin{equation*}
\begin{split}
\binom{k+1}{\ell+1}B_{n,k,\ell}&= \binom{k+1}{\ell+1}\tb_{k-1}(X)-\lambda_{\ell,k-\ell-2}(X) \\
&\leq \binom{k+1}{\ell+1}\tb_{k-1}(X-\sigma)-\lambda_{\ell,k-\ell-2}(X-\sigma)  \\
&=\binom{k+1}{\ell+1}B_{n,k,\ell}-\lambda_{\ell,k-\ell-2}(X-\sigma).
\end{split}
\end{equation*}
Therefore
\begin{equation}
\label{bxms}
\lambda_{\ell,k-\ell-2}(X-\sigma)=0.
\end{equation}
Combining \eqref{bkmo}, \eqref{bxms} and the minimality of $X$, it follows that
$\tb_k(X-\sigma)=0$. Using \eqref{lexact} for the complex $X-\sigma$ it follows that
\begin{equation*}
\begin{split}
f_k(X-\sigma)&= \tb_k(X-\sigma)-\tb_{k-1}(X-\sigma)+\binom{n-1}{k} \\
&=\binom{n-1}{k}-B_{n,k,\ell}=F_{n,k,\ell}.
\end{split}
\end{equation*}
Let $\tau \in \dn(\ell)$.  Using \eqref{lexact} for $\lk(X,\tau)$ we obtain
\begin{equation}
\label{eplk}
\tb_{k-\ell-1}(\lk(X-\sigma,\tau)- \tb_{k-\ell-2}(\lk(X-\sigma,\tau))=f_{k-\ell-1}(\lk(X-\sigma,\tau))-\binom{n-\ell-2}{k-\ell-1}.
\end{equation}
Summing \eqref{eplk} over all $\tau \in \dn(\ell)$ it follows that
\begin{equation}
\label{llklo}
\begin{split}
\lambda_{\ell,k-\ell-1}(X-\sigma)&=\lambda_{\ell,k-\ell-1}(X-\sigma)-\lambda_{\ell,k-\ell-2}(X-\sigma) \\
&=\sum_{\tau \in \dn(\ell)} f_{k-\ell-1}(\lk(X-\sigma,\tau))-\binom{n}{\ell+1}\binom{n-\ell-2}{k-\ell-1} \\
&=\binom{k+1}{\ell+1} f_k(X-\sigma)-\binom{n}{\ell+1}\binom{n-\ell-2}{k-\ell-1} \\
&=\binom{k+1}{\ell+1}F_{n,k,\ell}-\binom{n}{\ell+1}\binom{n-\ell-2}{k-\ell-1}=0.
\end{split}
\end{equation}
Choose a $k$-simplex $\sigma \neq \sigma' \in \supp(z)$ and an $\ell$-simplex $\tau \in \sigma'(\ell)-\sigma(\ell)$. Then on one hand
$\lk(X-\sigma,\tau)=\lk(X,\tau)$, hence by Claim \ref{ktol}
$$
\tb_{k-\ell-1}(\lk(X-\sigma,\tau))=\tb_{k-\ell-1}(\lk(X,\tau)) \neq 0.
$$
On the other hand it follows from \eqref{llklo} that $\tb_{k-\ell-1}(\lk(X-\sigma,\tau))=0$, a contradiction.
Thus $\tb_k(X)=0$ and hence
$$
f_k(X)= \tb_k(X)-\tb_{k-1}(X)+\binom{n-1}{k}=\binom{n-1}{k}-B_{n,k,\ell}=F_{n,k,\ell}.
$$
\noindent
(b) $\Rightarrow$ (c): Suppose $f_k(X)=F_{n,k,\ell}$. As $\lambda_{\ell,k-\ell-2}(X)=0$, it follows that
$\tb_{k-\ell-2}(\lk(X,\tau))=0$ and hence $f_{k-\ell-1}(\lk(X,\tau) \geq \binom{n-\ell-2}{k-\ell-1}$ for all $\tau \in \dn(\ell)$.
Therefore
\begin{equation}
\label{getff}
\begin{split}
F_{n,k,\ell}&=f_k(X)=\frac{1}{\binom{k+1}{\ell+1}}\sum_{\tau \in \dn(\ell)}f_{k-\ell-1}(\lk(X,\tau)) \\
&\geq \frac{\binom{n}{\ell+1}\binom{n-\ell-2}{k-\ell-1}}{\binom{k+1}{\ell+1}}=F_{n,k,\ell}.
\end{split}
\end{equation}
Hence $f_{k-\ell-1}(\lk(X,\tau))= \binom{n-\ell-2}{k-\ell-1}$ and therefore $\lk(X,\tau)$ is a $(k-\ell-1)$-hypertree for all $\tau \in \dn(\ell)$.
\\
(c) $\Rightarrow$ (a):  Assume that $\lk(X,\tau)$ is a $(k-\ell-1)$-hypertree for all $\tau \in \dn(\ell)$. Then, as in \eqref{getff}, it follows that
$f_k(X)=F_{n,k,\ell}$. Furthermore, by \eqref{bbkmo}
$$
\tb_k(X) \leq \frac{1}{\binom{k+1}{\ell+1}} \sum_{\tau \in X(\ell)} \tb_{k-\ell-1}(\lk(X,\tau))=0.
$$
Hence
$$
\tb_{k-1}(X)=\tb_k(X)-f_k(X)+\binom{n-1}{k}=\binom{n-1}{k}-F_{n,k,\ell}=B_{n,k,\ell}.
$$
{\enp}

\section{Links of Sum Complexes}
\label{s:sum}
Let $n$ be a prime and let $A$ be a subset of the cyclic group  $V=\Int_n$.
Identify the vertex set of $\dn$ with the elements of $\Int_n$. For $s \leq n-2$ define the {\it sum complex} $Y_{A,s+1} \subset \dn^{(s)}$
by
$$
Y_{A,s+1}=\dn^{(s-1)} \cup \{\sigma \subset \Int_n: |\sigma|=s+1 \text{~and~}\sum_{x \in \sigma} x \in A\}.
$$
The homology groups $\thh_*(Y_{A,s+1};\KK)$ were determined in \cite{LMR10,M14}.
When $A$ is a cyclic interval in $\Int_n$, the Betti numbers $\tb_*(Y_{A,s+1};\KK)$ do not depend on $\KK$ and take the following simple form.
\begin{theorem}[\cite{LMR10,M14}]
\label{sumcom}
Let $n$ be a prime and let $A=\{t,t+1,\ldots,t+r\}$ be an interval of size $r+1$ in $\Int_n$. Then for any field $\KK$
\begin{equation*}
\label{sum1}
\tb_i(Y_{A,s+1};\KK)=
\begin{cases}
\frac{s-r}{s+1}\binom{n-1}{s} & \text{if  } i=s-1,~ r \leq s,  \\
\frac{r-s}{s+1}\binom{n-1}{s} & \text{if  } i=s,~ r \geq s, \\
0 & \text{otherwise}.
\end{cases}
\end{equation*}
\end{theorem}
\noindent
Let $\ell \leq k-2$ and let $c_{k,\ell}=\frac{(\ell+1)(k-\ell)}{(k-\ell-1)!}$.
\begin{proposition}
\label{linkb}
Let $B=\{0,\ldots,k-\ell-1\}$. Then for any $\tau \in \dn(\ell)$
\begin{equation*}
\label{sum2}
\tb_{k-\ell-2}(\lk(Y_{B,k+1},\tau)) \leq c_{k,\ell} n^{k-\ell-2}.
\end{equation*}
\end{proposition}
\noindent
{\bf Proof:} Let $y=\sum_{x \in \tau} x$ and let $C=\{b-y:b \in B\}$. Then
$$
\lk(Y_{B,k+1},\tau)=Y_{C,k-\ell}[V-\tau] \subset Y_{C,k-\ell}.
$$
Applying Theorem \ref{sumcom} with $A=C$, $r=s=k-\ell-1$ and $i \in \{s-1,s\}$, it follows that
$$\tb_{k-\ell-2}(Y_{C,k-\ell})=\tb_{k-\ell-1}(Y_{C,k-\ell})=0.$$
Hence, the exact sequence
\begin{footnotesize}
\begin{equation*}
0=\thh_{k-\ell-1}(Y_{C,k-\ell})\rightarrow \thh_{k-\ell-1}(Y_{C,k-\ell},Y_{C,k-\ell}[V-\tau]) \rightarrow
\thh_{k-\ell-2}(Y_{C,k-\ell}[V-\tau]) \rightarrow  \thh_{k-\ell-2}(Y_{C,k-\ell})=0
\end{equation*}
\end{footnotesize}
implies
\begin{equation}
\label{isomor}
\thh_{k-\ell-1}(Y_{C,k-\ell},Y_{C,k-\ell}[V-\tau]) \cong
\thh_{k-\ell-2}(Y_{C,k-\ell}[V-\tau]).
\end{equation}
For $a,c \in \Int_n$ let
$$\cf_{a,c}=\{\sigma \in \dn(k-\ell-1): a \in \sigma\, , \, \sum_{x \in \sigma}x =c\}.$$
Note that any $\eta \in \dn(k-\ell-2)$ is contained in at most one $\sigma \in \cf_{a,c}$.
Therefore
\begin{equation}
\label{nkl2}
\begin{split}
|\cf_{a,c}|&=\frac{1}{k-\ell-1}|\{(\eta,\sigma): a \in \eta \subset \sigma \in \cf_{a,c} \, , \, |\eta|=k-\ell-1\}| \\
&\leq \frac{1}{k-\ell-1}|\{\eta \in \dn(k-\ell-2): a \in \eta\}|  \\
&= \frac{1}{k-\ell-1} \binom{n-1}{k-\ell-2} \leq \frac{n^{k-\ell-2}}{(k-\ell-1)!}.
\end{split}
\end{equation}
Combining \eqref{isomor} and \eqref{nkl2} it follows that
\begin{equation*}
\label{kminl}
\begin{split}
\tb_{k-\ell-2}(\lk(Y_{B,k+1},\tau))&=\dim \thh_{k-\ell-2}(Y_{C,k-\ell}[V-\tau]) \\
&=\dim \thh_{k-\ell-1}(Y_{C,k-\ell},Y_{C,k-\ell}[V-\tau]) \\
&\leq f_{k-\ell-1}(Y_{C,k-\ell}) - f_{k-\ell-1}(Y_{C,k-\ell}[V-\tau]) \\
&= |\{\sigma \in \dn(k-\ell-1): \sigma \cap \tau \neq \emptyset \, , \, \sum_{x \in \sigma} x \in C\}| \\
&\leq \sum_{(a,c)\in \tau \times C} |\cf_{a,c}|  \\
&\leq \frac{(\ell+1)(k-\ell)}{(k-\ell-1)!} n^{k-\ell-2}.
\end{split}
\end{equation*}
{\enp}

\section{The Lower Bound}
\label{s:lower}
\noindent
{\bf Proof of Theorem \ref{exacom}:} For the case $\ell=k-1$ see the remark in Section \ref{s:conc}.
Assume that $\ell \leq k-2$.
Let $n$ be a prime and let $B=\{0,\ldots,k-\ell-1\}$.
For any $\tau \in \dn(\ell)$
choose a set of $(k-\ell-1)$-dimensional simplices
$S_{\tau} \subset \lk(\dn,\tau)(k-\ell-1)$
of size
$|S_{\tau}|=\tb_{k-\ell-2}(\lk(Y_{B,k+1},\tau))$
such that
$$
\tb_{k-\ell-2}(\lk(Y_{B,k+1},\tau) \cup S_{\tau})=0.
$$
Proposition \ref{linkb} implies that
$$
|S_{\tau}| \leq c_{k,\ell}n^{k-\ell-2}.
$$
Let
$$
X_{n,k,\ell}=Y_{B,k+1} \bigcup \{\tau \cup\eta: \tau \in \dn(\ell), \eta \in S_{\tau}\}.
$$
Then for all $\tau \in \dn(\ell)$
$$
\tb_{k-\ell-2}(\lk(X_{n,k,\ell},\tau))=\tb_{k-\ell-2}(\lk(Y_{B,k+1},\tau) \cup S_{\tau})=0.
$$
Applying Theorem \ref{sumcom} with $A=B$, $r=k-\ell-1$, $s=k$ and $i=s-1$, it follows that
$$\tb_{k-1}(Y_{B,k+1})=\frac{\ell+1}{k+1} \binom{n-1}{k}.$$
Hence
\begin{equation*}
\label{sum3}
\begin{split}
\tb_{k-1}(X_{n,k,\ell}) &\geq \tb_{k-1}(Y_{B,k+1}) -\sum_{\tau \in \dn(\ell)} |S_{\tau}| \\
&\geq \frac{\ell+1}{k+1} \binom{n-1}{k} -\binom{n}{\ell+1}  c_{k,\ell} n^{k-\ell-2} \\
&\geq  \frac{\binom{n-1}{\ell} \binom{n-\ell-2}{k-\ell}}{\binom{k+1}{\ell+1}}(1-O(n^{-1})).
\end{split}
\end{equation*}
{\enp}

\section{Constructions of $J_{n,k}$ for $k \leq 3$}
\label{s:twoc}
In this section we describe the constructions that establish Theorem \ref{exactc}.
\ \\ \\
(i) Let $k=1$ and $n$ be even. Let $J_{n,1}$ be a perfect matching on the vertex set $[n]$.
Then $\tb_{-1}(\lk(J_{n,1},v))=0$ for all $v \in [n]$ and $\tb_0(J_{n,1})=\frac{n}{2}-1=B_{n,1,0}$.
\ \\ \\
(ii) Let $k=2$ and $n=3t+2$. Let $J_{n,2}$ be the $2$-dimensional complex on the vertex set $\Int_n$
(see Figure \ref{figure1}(a)) given by
$$
J_{n,2}=\dn^{(1)} \cup \left\{\{i,i+3j+1,i+3j+2\}: i \in \Int_n~,~ 0 \leq j \leq t-1 \right\}.
$$
\begin{proposition}
\label{ktwo}
$J_{n,2}$ satisfies $\tb_0(\lk(J_{n,2},v))=0$ for all $v \in \Int_n$
and $\tb_1(J_{n,2})=\frac{1}{3} \binom{n-2}{2}=B_{n,2,0}$.
\end{proposition}

\begin{figure}[h]
\subfloat[$J_{n,2}$]
{\label{fig:TN}
\scalebox{0.3}{\input{ttn1.pstex_t}}}
\hspace{50pt}
\subfloat[$\lk(J_{n,2},0)$]
{\label{fig:LTN}
\scalebox{0.4}{\input{link1.pstex_t}}}
%\vspace{-1cm}
\caption{}
\label{figure1}
\end{figure}
\begin{figure}[h]
\begin{center}
{\label{fig:TN}
\scalebox{0.3}{\input{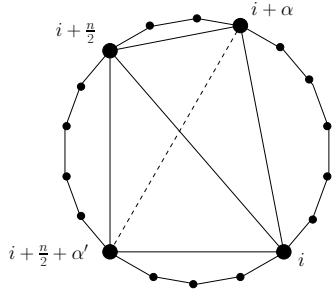}}}
\caption{A $3$-simplex in $J_{n,3}$}
\label{figure2}
\end{center}
\end{figure}

\noindent
{\bf Proof:}
By Theorem \ref{zlink2} it suffices to show that for all $i \in \Int_n$ the graph $\lk(J_{n,2},i)$ is a tree on the vertex set $\Int_n-\{i\}$.
By homogeneity it suffices to consider $\lk(J_{n,2},0)$.  It follows from the definition of $J_{n,2}$ that
$$\lk(J_{n,2},0)(1)= A_0 \cup B_0 \cup C_0$$
where
\begin{equation*}
\begin{split}
&A_0=\left\{\{3j+1,3j+2\}: 0 \leq j \leq t-1 \right\}, \\
&B_0=\left\{\{n-3j-1,1\}: 0 \leq j \leq t-1 \right\}=\left\{\{1,3j+1\}: 1 \leq j \leq t \right\},\\
&C_0=\left\{\{n-3j-2,n-1\}: 0 \leq j \leq t-1 \right\}=\left\{\{3j,n-1\}: 1 \leq j \leq t \right\}.
\end{split}
\end{equation*}
Thus, $\lk(J_{n,2},0)$ is the tree on $\Int_n-\{0\}$ depicted in Figure \ref{figure1}(b).
{\enp}
\ \\ \\
(iii) Let $k=3$ and $4 \leq n$ be even. Let $J_{n,3}$ be the $3$-dimensional complex
on the vertex set $\Int_n$ (see Figure \ref{figure2}) given by
$$
J_{n,3} = \dn^{(2)} \cup \left\{\{i,i+\alpha,i+n/2,i+\frac{n}{2}+\alpha'\}: 0 \leq i < \frac{n}{2}~,~ 0<\alpha,\alpha'<\frac{n}{2}\right\}.
$$
\begin{proposition}
\label{kthree}
$J_{n,3}$ satisfies $\tb_1(\lk(J_{n,3},v))=0$ for all $v \in \Int_n$
and $\tb_2(J_{n,3})=\frac{1}{4} \binom{n-2}{3}=B_{n,3,0}$.
\end{proposition}

\begin{figure}[h]
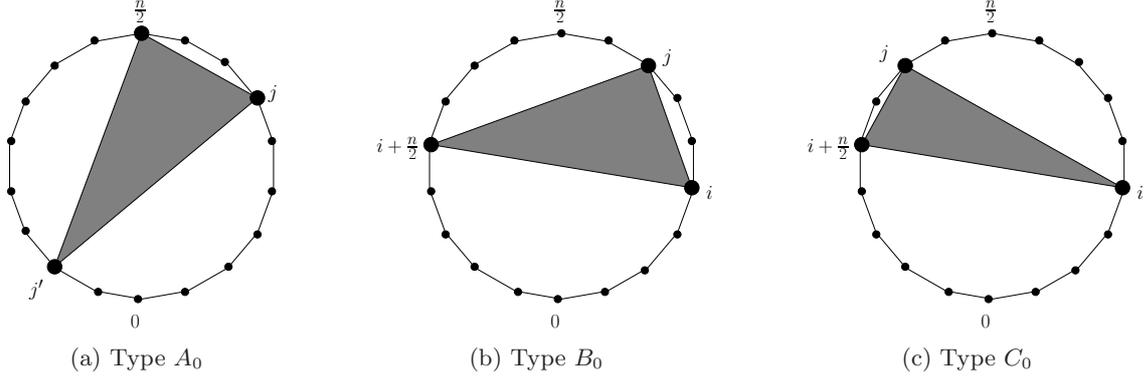

\subfloat[Type $A_0$]
{\label{fig:TNa}
\scalebox{0.3}{\input{jjna.pstex_t}}}
\hspace{30pt}
\subfloat[Type $B_0$]
{\label{fig:LTb}
\scalebox{0.3}{\input{jjnb.pstex_t}}}
\hspace{30pt}
\subfloat[Type $C_0$]
{\label{fig:LTc}
\scalebox{0.3}{\input{jjnc.pstex_t}}}
\caption{Three types of $2$-simplices in $\lk(J_{n,3},0)$}
\label{figure3}
\end{figure}

\noindent
{\bf Proof:} As in the proof of Proposition \ref{ktwo}, it suffices to show that $\lk(J_{n,3},0)$ is a $2$-hypertree on the vertex set $\Int_n-\{0\}$. We claim that $\lk(J_{n,3},0)$ is in fact collapsible. Partition the $2$-simplices of $\lk(J_{n,3},0)$ into $3$ disjoint families
(see Figure \ref{figure3}):
$$\lk(J_{n,3},0)(2)=A_0 \cup B_0 \cup C_0,$$
where
\begin{equation*}
\begin{split}
&A_0=\left\{ \{ j,n/2,j' \} : 0<j<n/2<j'<n \right\}, \\
&B_0=\left\{ \{ i,j,i+n/2 \} : 0<i<j<n/2 \right\}, \\
&C_0=\left\{ \{ i,j,i+n/2 \} : 0<i<n/2<j<i+n/2 \right\}.
\end{split}
\end{equation*}
If $0<i<j<n/2$ then the edge $\{i,j\}$ is contained in the unique $2$-simplex $\{i,j,i+n/2\} \in B_0$.
If $0<i<n/2<j<i+n/2$ then the edge $\{j,i+n/2\}$ is contained in the unique $2$-simplex $\{i,j,i+n/2\} \in C_0$.
Collapsing all these edges and the corresponding $2$-simplices, the resulting complex consists of all simplices in $A_0$ and their faces.
This complex is a cone on the vertex $n/2$ and is therefore collapsible.
{\enp}

\section{Concluding Remarks}
\label{s:conc}

We have shown that if $\dn^{(k-1)} \subset X \subset \dn^{(k)}$ satisfies
$\tb_{k-\ell-2}(\lk(X,\tau);\KK)=0$ for all $\tau \in \dn(\ell)$, then $\tb_{k-1}(X;\KK) \leq B_{n,k,\ell}$.
Furthermore, this bound is asymptotically tight for fixed $k,\ell$ and $n \rightarrow \infty$, and exact for $(k,\ell)=(1,0),(2,0),(3,0)$ and infinitely many $n$'s.
We suggest the following
\begin{conjecture}
\label{exx}
For any fixed $0 \leq \ell <k$ there exists a constant $n_0(k,\ell)$ such that if $n \geq n_0(k,\ell)$ and if $B_{n,k,\ell}$ is an integer, then
there exists a complex $\dn^{(k-1)} \subset X \subset \dn^{(k)}$ such that $\lambda_{\ell,k-\ell-2}(X)=0$ and $\tb_{k-1}(X)=B_{n,k,\ell}$.
\end{conjecture}
\noindent
{\bf Remarks:} 
\\
1. Let $\ell=k-1$. It follows from Theorem \ref{zlink2} that $\dn^{(k-1)} \subset X \subset \dn^{(k)}$ satisfies both
$\lambda_{\ell,k-\ell-2}(X)=\lambda_{k-1,-1}(X)=0$ and $\tb_{k-1}(X)=B_{n,k,k-1}=(\frac{n}{k+1}-1)\binom{n-1}{k-1}$ iff each $(k-1)$-simplex $\tau \in \dn(k-1)$
is contained in a unique $k$-simplex in $X(k)$, i.e. iff $S=X(k) \subset \binom{[n]}{k+1}$ is a Steiner system of type $S(k,k+1,n)$. Thus the case $\ell=k-1$ of Conjecture \ref{exx}
follows from Keevash groundbreaking work \cite{Keevash14} on the existence of designs.
\\
2. It would be interesting and useful for various applications to interpolate between Corollary \ref{bzlink} and Garland's Theorem, in particular to obtain sharp upper bounds on the rational Betti number $\tb_{k-1}(X;\Rat)$ in terms of $\tilde{\mu}_{k-\ell-2}(X)=\min_{\tau \in X(\ell)} \mu_{k-\ell-2}(\lk(X,\tau))$
when $0< \tilde{\mu}_{k-\ell-2}(X) \leq \frac{\ell+1}{k}$.
\\
3. The characterization given in Theorem \ref{zlink2} and the examples in Section \ref{s:twoc} suggest some natural questions concerning hypertrees, e.g. for which $(k-1)$-hypertrees
$\Delta_{n-2}^{(k-2)} \subset T \subset \Delta_{n-2}^{(k-1)}$ there exists a complex $\dn^{(k-1)} \subset X \subset \dn^{(k)}$ such that
$\lk(X,v) \cong T$ for all vertices $v \in X(0)$ ?

\end{document}